\documentclass[lettersize,joof thenal]{IEEEtran}
\usepackage{amsmath,amsfonts}
\usepackage{algorithmic}
\usepackage{titlesec}
\usepackage{algorithm}
\usepackage{soul}
\usepackage{mathtools}
\usepackage{subfigure}
\usepackage{amsthm}
\usepackage{amssymb}
\usepackage{amsmath}
\usepackage{array}
\usepackage{bm}
\usepackage[caption=false,font=normalsize,labelfont=sf,textfont=sf]{subfig}
\usepackage{textcomp}
\usepackage{stfloats}
\usepackage{url}
\usepackage{verbatim}
\usepackage{graphicx}

\usepackage{cite}
\usepackage[T1]{fontenc}
\usepackage{soul}
\usepackage{xcolor}

\begin{document}

% \title{A Sample Article Using IEEEtran.cls\\ for IEEE Journals and Transactions}

%Connections
\title{Slow, Fast and Opportunistic FAMA: A Spatial Block-Correlation Analysis under Nakagami-$m$ Fading Channels}

        % <-this % stops a space

\author{ 
Paulo R. de Moura,~\IEEEmembership{Graduate Student Member,~IEEE,} 
Hugerles S. Silva,~\IEEEmembership{Senior Member,~IEEE},\\
Ugo S. Dias,~\IEEEmembership{Senior Member,~IEEE}, 
Higo T. P. Silva,~\IEEEmembership{Member,~IEEE}
% Rausley A. A. de Souza,~\IEEEmembership{Senior Member,~IEEE} 
% and Osamah Badarneh,~\IEEEmembership{Senior Member,~IEEE}
 \thanks{P. R. de Moura is with the Electrical Engineering Department, University of Brasília, Brasília 70910-900, Brazil, and also with the Technical Advisory Unit, National Telecommunications Agency (Anatel), Brasília 70070-940, Brazil (e-mail: paulomoura@anatel.gov.br).}
%  \thanks{H. S. Silva is with the Electrical Engineering Department, University of Brasília (UnB), Federal District, Brazil, and also with Instituto de Telecomunicações and Departamento de Eletrónica, Telecomunicações e Informática, Universidade de Aveiro, Campus Universitário de Santiago, 3810-193 Aveiro, Portugal. (e-mail: hugerles.silva@av.it.pt)}
%  \thanks{U. S. Dias is with the Electrical Engineering Department, University of Brasília (UnB), Federal District, Brazil. (e-mail: udias@unb.br)}
% \thanks{H. T. P. Silva is with the Electrical Engineering Department, University of Brasília (UnB), Federal District, Brazil. (e-mail: higo.silva@unb.br)}
\thanks{H. S. Silva, U. S. Dias and H. T. P. Silva are with the Electrical Engineering Department, University of Brasília (UnB), Brasília 70910-900, Brazil.}
}

% \thanks{R. A. A. de Souza is with the National Institute of Telecommunications (Inatel), Santa Rita do Sapucaí 37540-000, Brazil. (e-mail: rausley@inatel.br)}
% \thanks{O. Badarneh is with the Electrical Engineering Department, School of Electrical Engineering and Information Technology, German Jordanian University, Amman 11180, Jordan. (e-mail: osamah.badarneh@gju.edu.jo)}

% }

% \markboth{Submitted to IEEE, April~2025}%
% {Moura \MakeLowercase{\textit{et al.}}:  Opportunistic fast and slow FAMA under Nakagami-$m$ Fading Channels}

\maketitle

\begin{abstract}
This paper studies slow, fast and opportunistic fluid antenna multiple access (FAMA) under the effect of Nakagami-$m$ fading channels, considering the new and realistic spatial block-correlation model.
Expressions for the outage probability (OP), based on the signal-to-interference ratio (SIR), are derived for slow FAMA. Interestingly, we provide
 mathematical relationships that allow the expressions of fast FAMA to be obtained from slow FAMA.
 %of one type of FAMA to be obtained from the other.
Multiplexing gains for an opportunistic FAMA (O-FAMA) network are presented
for both slow and fast FAMA scenarios.
Our analytical results are validated through Monte Carlo simulations,  under various channel and system parameters. 
All expressions derived in this work are original.
%and have not been previously published. 

\end{abstract}

\begin{IEEEkeywords}
FAMA,  multiplexing gain, Nakagami-$m$ fading, opportunistic FAMA, spatial block-correlation.
\end{IEEEkeywords}

\section{Introduction}

\IEEEPARstart{F}{luid} antenna systems (FAS) is a disruptive technology that promises to help to achieve the demand for massive connectivity of emerging mobile communication systems~\cite{Wong2020_BruceLee, Wong2021FAS}. 
A fluid antenna (FA) is a flexible, electronically reconfigurable antenna structure based on liquid or pixel technology.
In its canonical form, it consists of a linear structure with predefined positions, known as ports, where the radiating element is switched to optimize a reception metric. The performance benefits of FAS have been extensively analyzed in various operational modes and channel fading scenarios~\cite{Xu2024_Secret_FAS,Ghadi2024_Cache_FAS,Wong2021FAS}.
%\footnote{Interested readers may refer to \cite{New2024_FAS6G_Tutorial} for a comprehensive tutorial on FAS for 6G networks, which summarizes recent researches.}. 

Extending the concept of FAS, 
fluid antenna multiple access (FAMA) introduces multiple users equipped with FAs~\cite{Wong2022FAMA, Wong2023FastFamaMassiveConnectivity, Wong2023SlowFAMA} sharing the same resources.
In an FAMA scheme, FAs dynamically reconfigure their ports to mitigate interference in a shared spectrum scenario and maximize the signal-to-interference ratio (SIR) or signal-to-interference-plus-noise ratio (SINR). 
By leveraging the fading depth across the FA space, user equipments (UEs) can select ports that enhance performance in interference-limited environments. 
The literature identifies two primary FAMA operating modes: fast FAMA ($f$-FAMA)~\cite{Wong2022FAMA, Wong2023FastFamaMassiveConnectivity} and slow FAMA ($s$-FAMA)~\cite{Wong2023SlowFAMA}. 
While $f$-FAMA provides substantial performance gains, its practical applicability is constrained by the requirement to switch FAs at symbol time~\cite{Wong2023FastFamaMassiveConnectivity}. 
In turn, $s$-FAMA offers a more feasible alternative, requiring switching only between channel coherence times, when significant channel variations occur~\cite{Wong2023SlowFAMA}. 

Recently, another proposed approach for FAMA is based on opportunistic scheduling. 
This technique relies on the dynamic allocation of resources to a subset of UEs from a large user pool based on their channel conditions~\cite{Wong2023OpportunisticFAMA}. 
The so-called opportunistic-FAMA (O-FAMA) integrates opportunistic scheduling with either $f$-FAMA or $s$-FAMA, leveraging their respective advantages to improve network performance~\cite{Wong2023OpportunisticFAMA}.
Using a reinforcement learning approach, it is shown in~\cite{Waqar2024_Opportunistic_FAMA_RL} that it is possible to select the best users and FAS ports to reach a network sum rate close to the ideal, which makes O-FAMA an interesting option for multiple access.

In the literature,~\cite{Wong2022FAMA,Wong2023SlowFAMA,Wong2023OpportunisticFAMA,Waqar2024_Opportunistic_FAMA_RL} assume Rayleigh fading channels and \cite{Wong2023FastFamaMassiveConnectivity} considers the finite-scatterer channel model.
Despite serving as a fundamental benchmark, the Rayleigh distribution offer no additional degrees of freedom~(DoF), significantly limiting the applicability of their results. 
More recently, only~\cite{Moura2025_sFAMA_Nakagami} extended the analysis of $s$-FAMA to Nakagami-$m$ fading.
Furthermore, the correlation models originally used in the analysis of~\cite{Wong2022FAMA,Wong2023SlowFAMA,Wong2023OpportunisticFAMA} are based on restricted forms of the correlation matrix \cite{Beaulieu2011, Wong2022Closed} for easy mathematical tractability. However, they do not capture well the physical behavior of FAS. 
New correlated channel models are also presented in~\cite{Khammassi2023} and~\cite{New2024_FAS_New_Insights}, that have been proven to be prohibitively complex, resulting in intractable analysis in FAS and FAMA \cite{Khammassi2023, New2024_FAS_New_Insights, Xu2023TwoUserFAMA}, due to multi-folded integrals involved.

In this context, we study for the first time the performance of $s$-FAMA, $f$-FAMA and O-FAMA under Nakagami-$m$ fading channels considering a spatial block-correlation analysis.
We adopted the Nakagami-$m$ fading since it is a well-established and relevant model for evaluating the performance of current and emerging systems, providing greater flexibility compared to Rayleigh fading.
We also adopted the spatial block-correlation model recently proposed in~\cite{RamirezEspinosa2024_BlockCorrelationFAS}, which accurately characterizes the correlation, as predicted by classical realistic models such as Jakes’s, while maintaining analytical tractability and simplicity of the constant correlation model presented in~\cite{Wong2022Closed}.
To the best of our knowledge, the analyses and all the expressions presented here are novel in the literature.

The main contributions of this article are summarized as:
\begin{itemize}
     \item Novel and more precise expressions are derived for outage probability~(OP)-based on SIR for $s$-FAMA under the effect of Nakagami-$m$ fading, considering the spatial block-correlation model.
     \item A new and approximate expression is derived for the SIR-based OP under $s$-FAMA, over Nakagami-$m$ channels with spatial block-correlation. Furthermore, an upper bound for the SIR-based OP is also deduced.
     \item Interesting similarities between fast and slow FAMA are identified, where mathematical relationships are presented. These relationships allow expressions of one type of FAMA to be obtained from the other.
     \item O-FAMA is analyzed for fast and slow FAMA, under Nakagami-$m$ channels and spatial block-correlation, where results are provided for the multiplexing gain.
\end{itemize}

%The remainder of this paper is organized as follows. The analysis for fast, slow and opportunistic FAMA under Nakagami-$m$ channels considering spatial block-correlation model are shown in Sections~\ref{sFAMA},~\ref{fFAMA} and~\ref{Ofama} respectively. Section~\ref{results} shows the numerical results. Section \ref{conclusions} brings the paper's conclusions.

\section{FAMA Model} \label{sec:system_model}

\subsection{System Model}

\begin{figure}[!t]
    \centering
    \includegraphics[width=1\linewidth]{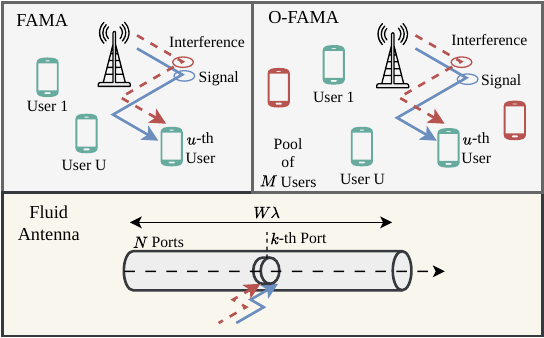}
    \caption{FAMA and O-FAMA systems.
    }
    \label{fig:fama_system}
\end{figure}

The downlink FAMA network considered in this work is illustrated in Fig.~\ref{fig:fama_system} and consists of a BS equipped with $U$ antennas. In this system, each antenna of  BS transmits a signal destined to a specific user of the network, thus constituting a network with $U$ users, and each UE contains an FA with $N$ ports.
The received signal at the $n$-th port, of a FAS for a given user $u$, is modeled as
\begin{align}\label{eq:model_Fama}
     r^{(u)}_{n} = s_{u} h^{(u,u)}_{n}+  \sum _{\substack{ \tilde{u}\ne u}}^{U} s_{\tilde {u}}  h^{(\tilde{u},u)}_{n} + \eta _{n}^{(u)},
\end{align}
in which $s_u$ denotes the transmitted symbol intended for the $u$-th user, $h^{(u,u)}_{n}$ is the corresponding complex fading channel experienced at $n$-th port of  user $u$, and  $h^{(\tilde{u},u)}_{n}$ denotes the fading channel from the BS antenna transmitting user $\tilde{u}$’s signal, $s_{\tilde u}$, which acts as an interference at $n$-th port of user $u$.
Furthermore,  $\eta^{(u)}_n$ is the complex additive white Gaussian noise (AWGN), at the $n$-th port for user $u$, with zero mean and variance $\sigma^2_{\eta}$. 
Note that each BS antenna is assigned to transmit the signal for a given user on the downlink. 
The average power of the transmitted symbol is $\sigma_s^2 = \mathbb{E}[|s_u|^2]$, $\forall u$, in which $\mathbb{E}[\cdot]$ is the expectation operator. In this paper, we assume that the UEs have perfect knowledge of the channel, so the best ports are selected; and also that the switching delay between ports is negligible.

\subsection{Channel and Spatial Correlation Models}

The channel envelope for the one-dimensional (1D) FAs under Nakagami-$m$ fading is expressed as
\begin{equation} \label{eq:hn}
    \lvert h^{(u,u)}_{n} \rvert = \sqrt{\sum_{l=1}^m {\lvert g^{(u,u)}_{nl} \rvert^2}},
\end{equation}
where $m$ is the fading severity and $\{g^{(u,u)}_{nl}\}$ is a set of mutually correlated complex Gaussian random variables (RVs) with zero mean and $\mathbb{E}(\lvert g^{(u,u)}_{nl} \rvert^2) = 2\sigma_u$. In general, the mutual correlations between the channel coefficients of any pair of ports are described by the spatial correlation matrix $\bm{\Sigma} \in \mathbb{C}^{N \times N}$. For one-dimensional FAs in an isotropic scattering environment, under the Jakes model, the elements of the correlation matrix are given by
$    [\bm{\Sigma}]_{nk} = \text{J}_0 \left( \frac{2\pi (n-k) W}{N-1} \right)$,
% \begin{equation}
%     [\bm{\Sigma}]_{nk} = \text{J}_0 \left( \frac{2\pi (n-k) W}{N-1} \right),
% \end{equation}
where $W$ denotes the normalized antenna size and $\text{J}_0(\cdot)$ is the zeroth-order Bessel function of the first kind. Alternative spatial correlation models have been proposed in the literature to account for various antenna geometries and multipath scattering distributions~\cite{Queiroz2011_Spatial_Correlation_DoA}. While these models effectively capture the mutual correlation properties between ports, the analytical tractability of the FAS and FAMA systems performance is prohibitively complex~\cite{Khammassi2023, New2024_FAS_New_Insights, Xu2023TwoUserFAMA}. %\cite{RamirezEspinosa2024_BlockCorrelationFAS}.

To mitigate this complexity while preserving the correlation effects along the FAs, we adopt the methodology proposed in~\cite{RamirezEspinosa2024_BlockCorrelationFAS}, which approximates the correlation matrix by a block-diagonal matrix with constant correlation coefficients $\delta$, denoted by $\hat{\bm{\Sigma}}$. In this approach, the $N$ ports of a FA are partitioned into $B$ blocks, where each block has length of $L_b$ ports, with $b \in [1, \cdots, B]$ and $\sum_{b=1}^B L_b = N$. According to the algorithm in~\cite{RamirezEspinosa2024_BlockCorrelationFAS}, both $L_b$ and $\hat{\bm{\Sigma}}$ are determined from the dominant eigenvalues of the reference spatial correlation matrix ${\bm{\Sigma}}$ with respect to a threshold $\rho_{\text{th}}$. Consequently, $\hat{\bm{\Sigma}}$ is constructed as a block-diagonal matrix consisting of $B$ equi-correlation submatrices of size $L_b \times L_b$. Based on this representation, the OP of the FAMA system can be evaluated as $B$ independent blocks, while approximately retaining the correlation properties. There is no closed-form solution for the best choice of $\delta$ and $\rho_{\text{th}}$ to optimize $\hat{\bm{\Sigma}}$, but in \cite{RamirezEspinosa2024_BlockCorrelationFAS} it is suggested to choose $\delta \in (0.95,0.99)$ and $\rho_{\text{th}}=1$. In principle, this method is applicable to any reference correlation matrix, providing flexibility across different scattering scenarios.

\section{Slow FAMA under Nakagami-$m$ Channels with Spatial Block-Correlation Model} \label{sFAMA}

\subsection{Channel Model}

Based on Section~\ref{sec:system_model}, the channel coefficients are defined as
\begin{align}\label{eq:gk_sFAMA}
 g^{(u,u)}_{nl}
 &= \sigma_u \left({\sqrt {1-\delta} \, x^{( u,u)}_{nl}+ \sqrt{\delta} \,
 \mathrm{x}^{(u,u)}_{b(n)l}}\right) \nonumber \\
& +j \sigma_u \left ({\sqrt{1-\delta} \, y^{(u,u)}_{nl}+ \sqrt{\delta} \,
\mathrm{y}^{( u,u)}_{b(n)l}}\right),
\end{align}
where $\mathrm{x}^{(u,u)}_{b(n)l}$, $x^{(u,u)}_{1l}$,$\dots$, $x^{(u,u)}_{Nl}$ and $\mathrm{y}^{(u,u)}_{b(n)l}$, $y^{(u,u)}_{1l}$,$\dots$, $y^{(u,u)}_{Nl}$, with $l \in [1,\dots,m]$, are zero mean and unit variance independent Gaussian RVs. The RVs $\mathrm{x}^{(u,u)}_{b(n)l}$ and $\mathrm{y}^{(u,u)}_{b(n)l}$ are referenced to the block index $b(n)$, which is defined as $b(n)=1$, for $n=1,\dots,L_1$, $b(n)=2$, for $n=L_1+1,\dots,L_2$ and so on. Its is assumed that $\sigma_{u} = \sigma, \, \forall \tilde u$.
In (\ref{eq:gk_sFAMA}), $\delta$ is the common power correlation coefficient between any two ports in a block.

Note that the expressions (\ref{eq:hn}) and (\ref{eq:gk_sFAMA}) are also suitable for interfering users. In this case, the distribution parameters are indicated by $m_{\tilde u}$, with $l \in \{ 1,\dots,m_{\tilde u} \} $, $\sigma_{\tilde u} = \sigma~\forall u$, and the superscripts are denoted by $(\tilde u, u)$.

\subsection{SIR Model}

Assuming an interference-dominated scenario, where noise can be neglected in~(\ref{eq:model_Fama}), a $s$-FAMA UE selects the port where the SIR is maximized, such as
\begin{equation} \label{eq:SIR_BC}
    \text{SIR}= \underset{n}{\text{max}} \frac{\sigma_s^2 \lvert h_n^{(u,u)} \rvert^2}
    {\sigma_s^2 \sum_{\tilde u \ne u}^U \lvert h_n^{(\tilde u,u)} \rvert^2}
    = \underset{n}{\text{max}} \frac{X_n}{Y_n},
\end{equation}
in which
\begin{equation} \label{eq:Xk}
    X_n  = \sum_{l=1}^{m} 
    \left(x^{(u,u)}_{nl}+ \varphi \,  \mathrm{x}^{(u,u)}_{b(n)l} \right)^2   
    + \left(y^{(u,u)}_{nl}+\varphi \,  \mathrm{y}^{(u,u)}_{b(n)l}\right)^2
\end{equation}
and
\begin{align} \label{eq:Yk}
    Y_n  =  
     \sum_{\tilde u \ne u}^{U} \sum_{l=1}^{m_{\tilde u}}  
     \left(x^{(\tilde u,u)}_{nl}+ \varphi \, \mathrm{x}^{(\tilde u,u)}_{b(n)l} \right)^2 
    + \left(y^{(\tilde u,u)}_{nl}+ {\varphi} \, \mathrm{y}^{(\tilde u,u)}_{b(n)l} \right)^2
\end{align}
with $\varphi = \sqrt{\delta/(1-\delta)}$. 
% Note that the two sums in $Y_n$ can be combined in one sum, such as $\sum_{\tilde u \ne u}^{U} \sum_{l=1}^{m_{\tilde u}} = \sum_{\tilde u \ne u}^{\widetilde U}$, where $\widetilde U = \sum_{\tilde u \ne u}^U m_{\tilde u}$.

\subsection{Outage Probability Analysis}

\subsubsection{Exact Expression}
The OP, considering that the blocks are independent, is given by
\begin{equation} \label{eq:OP_block}
    {P}_{\mathrm{out}} \triangleq \text{Pr}[\text{SIR} \leq \gamma] = \prod_{b=1}^B P_{\mathrm{out};b} (\gamma),
\end{equation}
in which $P_{\mathrm{out};b}(\gamma)$ is the OP of the $b$-th block and $\gamma$ is the SIR threshold.

Replacing \cite[Eq. (23)]{Moura2025_sFAMA_Nakagami} into (\ref{eq:OP_block}), the SIR-based OP results in
\begin{align} \label{eq:OP_SIR_sFAMA}
    {P}_{\text {out}}  \!=\!
    \prod _{b=1}^{B} \int _{0}^{\infty }\int _{0}^{\infty } &
    \frac {r_b^{m-1} \widetilde {r}_{b}^{\widetilde U-1}e^{-\frac {r_{b} \!+\! \widetilde {r}_{b}}{2}}}
    {2^{m+ \widetilde U}\Gamma(m) \Gamma (\widetilde U)} \nonumber \\
    &\times\left [{{G(\gamma ; r_{b}, \widetilde {r}_{b})}}\right ]^{L_{b}}\,\text {d} r_{b} \,\text {d}\widetilde {r}_{b}, 
\end{align}
where $G(\gamma ; r_{b}, \widetilde {r}_{b})$ is given by (\ref{eq:G_OP}), $(x)_j$ denotes the Pochhammer symbol, $\text{Q}_{\nu}(\cdot)$ is the $\nu$-th order Marcum-Q function and $\text{I}_{\nu}(\cdot)$ is the $\nu$-th order modified Bessel function of the first kind. 

\subsubsection{Approximation}
Applying the generalized Gauss-Laguerre quadrature \cite{Weisstein2025_Laguerre_Gauss_Quadrature} in~(\ref{eq:OP_SIR_sFAMA}), we have a simple approximation for the SIR-based OP, as 
\begin{align} \label{eq:OP_SIR_sFama_Quadrature}
    {P}_{\text {out}} \approx \prod _{b=1}^{B} 
    \tfrac {1}{\Gamma(m) \Gamma (\widetilde U)}
    \sum _{i=1}^{n_I} \sum _{j=1}^{n_J} w_i  w_j \left [{{G(\gamma ; 2\,x_i, 2\,x_j)}}\right ]^{L_{b}}, 
\end{align}
where $x_i$, for $i\in [1,\dots,n_I]$, and $x_j$, for $j\in[1,\dots,n_J]$, are, respectively, the roots of the generalized Laguerre polynomials $L^{m-1}_{n_I}(x_i)$ and $L^ {\tilde U-1}_{n_J} (x_j)$, with weights 
$w_i = \frac{ \Gamma(n_I+m) x_i}{ n_I!  (n_I+1)^2 [L^{m-1}_{{n_I}+1} (x_i)  ]^2 }$ 
and $w_j = \frac{ \Gamma(n_J + \tilde U) x_j}{ n_J!  (n_J+1)^2 [L^{\tilde U-1}_{{n_J}+1} (x_j) ]^2 }$.

It should be mentioned that the use of the quadrature technique reduces computational time while maintaining an acceptable accuracy in the OP calculation.

\subsubsection{Upper Bound}
The upper bound for the OP can be derived using the fact that for $\delta \rightarrow 1$, the exponential in (\ref{eq:G_OP}) tends to zero, and then the OP is approximated only by the term $\text{Q}_{\tilde U}(\cdot)$.
Furthermore, for very large $N$ that results in large $L_b$ for most dominant eigenvalues, it follows that $[\text{Q}_{\tilde U}(\cdot)]^{L_p}$ tends to a Heaviside step function, shift by a threshold $f(\widetilde r_b)$. Thus, $[G(\gamma ; r_{b}, \widetilde {r}_{b}) ]^{L_b} = 1$ for $r_b < f(\widetilde r_b)$ and $[G(\gamma ; r_{b}, \widetilde {r}_{b}) ]^{L_b} = 0$ for $r_b > f(\widetilde r_b)$, so (\ref{eq:OP_SIR_sFAMA}) is approximated by 
\begin{equation} \label{eq:OP_SIR_sFAMA_step_approx}
{P}_{\text {out}}  \approx
    \prod _{b=1}^{B} \int _{0}^{\infty } 
    \frac {\widetilde {r}_{b}^{\widetilde U-1}e^{-\widetilde r_b/2 }}
    {2^{m+ \widetilde U}\Gamma(m) \Gamma (\widetilde U)} 
    \int _{0}^{ f(\widetilde r_b) }  
    {r^{m-1} e^{-\frac {r_{b} }{2}}}
    \text {d} r_{b} \,\text {d}\widetilde {r}_{b} .
    \end{equation}
Setting the threshold as $f(\widetilde r_b) = \gamma \, \widetilde  r_b$ and using \cite[Eqs. (3.351-1) and (3.351-3)]{Gradshteyn2000}, the integrals in (\ref{eq:OP_SIR_sFAMA_step_approx}) are solved, which results in the upper bound
\begin{equation} \label{eq:OP_upper_bound}
    {P}^{\text{UB}}_{\text {out}} \approx  
    \left( 1 - \sum_{i=0}^{m-1} \frac{\gamma^i (\widetilde U)_i}{i!(\gamma+1)^{\widetilde U + i} } \right)^B.
\end{equation}
It can be shown that (\ref{eq:OP_upper_bound}) is also the OP of $B$ independent antennas under Nakagami-$m$ channels.

\begin{figure*}
\begin{align} \label{eq:G_OP}
    G(\gamma ; r_{b}, \widetilde {r}_{b}) 
    & =\; \text{Q}_{\widetilde U}\left ({{\sqrt {\frac {\delta \, \gamma \, \widetilde {r}_{b}}{(1-\delta)(\gamma +1)}}, 
    \sqrt {\frac {\delta \, r_{b}}{(1-\delta)(\gamma +1)}}}}\right)  - 
    \left ({{\frac {1}{\gamma +1}}}\right)^{m + \widetilde U-1} \left( \frac{r_b}{\gamma \, \widetilde r_b } \right)^{\frac{1-m}{2}}
    \exp \left ({{-\frac {\delta}{2(1-\delta)}\frac {\gamma \widetilde {r}_{b}+r_{b}}{\gamma +1}}}\right) \nonumber \\ 
    & \quad \times \,\sum _{k=0}^{m+\widetilde U-2}\sum _{j=0}^{m+ \widetilde U-k-2}\frac {(m+\widetilde U-(j+k)-1)_{j}}{j!}\left ({{\frac {r_{b}}{\widetilde {r}_{b}}}}\right)^{\frac {j+k}{2}} 
    (\gamma +1)^{k}\gamma ^{\frac {j-k}{2}}\text{I}_{1-m+j+k}
    \left ({{\frac {\delta \sqrt{\gamma r_{b}\widetilde {r}_{b}}}{(1-\delta)(\gamma +1)}}}\right) 
    \end{align}
    \hrulefill\vspace{-0.5cm}
\end{figure*}

\subsection{Multiplexing Gain}

The multiplexing gain for $s$-FAMA, denoted as $\mathcal{G}_m$, can be defined as \cite[Eq. (29)]{Wong2023SlowFAMA}
\begin{equation} \label{eq:mult_gain}
    \mathcal{G}_m = U\,(1-P_\mathrm{out}).
\end{equation}

\section{Fast FAMA under Nakagami-$m$ Channels with Spatial Block-Correlation Model} \label{fFAMA}

\subsection{Channel Model}

For $f$-FAMA, the total interference in the received signal is treated as a single RV $\widetilde h_n^{(u)}=\sum_{\tilde u \ne u}^U s_{\tilde u} h_n^{(\tilde u, u)}$. Considering the seminal work of Nakagami \cite{Nakagami1960}, the sum of the complex Nakagami-$m$ interference RVs is approximated by another Nakagami-$m$  RV, with average power $\widetilde \Omega$ and fading parameter $\widetilde m$ given by \cite[Eq. (96)]{Nakagami1960}
\begin{equation} \label{eq:OmegaTil}
    \widetilde \Omega =  \sum_{\tilde u \ne u}^U  \Omega_{\tilde u}  \sigma_s^2 = 
    2 \sigma^2   \sigma_s^2  \widetilde U
\end{equation}
and
\begin{align} \label{eq:mTil}
    \widetilde m = 
      \frac{\left(\sum_{\tilde u \ne u}^U   m_{\tilde u}\right)^2}
     {\sum_{\tilde u \ne u}^U m_{\tilde u} + \sum_{\tilde u \ne u}^U  \sum_{\substack{i\neq\tilde{u} \\ \tilde{u}\neq u}}^{U} m_{\tilde u} m_i } ,
\end{align}
in which $\Omega_{\tilde u}=\mathbb{E}[\lvert h^{(\tilde{u},u)}_n \rvert^2]= 2 \, m_{\tilde u} \, \sigma^2$ and $\mathbb{E}[s_{\tilde u}^2] = \sigma_s^2, \, \forall \tilde u$. 
Therefore, the total interference power is modeled as
$\lvert \widetilde h^{(u)}_{n}  \rvert^2 = \lvert \sum_{\tilde u \ne u}^U s_{\tilde u} h_n^{(\tilde u, u)} \rvert^2 \approx \sum_{l=1}^{\widetilde m} \lvert  \widetilde{g}^{(u)}_{n  l} \rvert ^2  $ 
with $\widetilde{g}^{(u)}_{n  l}$ being expressed similarly to (\ref{eq:gk_sFAMA}), with Gaussian components denoted as $\widetilde {\mathrm{x}}^{(u)}_{b(n)l}, \widetilde x^{(u)}_{1l},\dots, \widetilde x^{(u)}_{Nl}$ and $\widetilde {\mathrm{y}}^{(u)}_{b(n)l}, \widetilde y^{(u)}_{1l},\dots, \widetilde y^{(u)}_{Nl}$ and $\mathbb{E}[\rvert\widetilde{g}^{(u)}_{n l}\lvert^2] = 2\widetilde\sigma$ and $\widetilde \Omega = \mathbb{E}[\lvert \widetilde{h}^{(u)}_n \rvert^2] = 2  \widetilde{m}  \widetilde{\sigma}^2$.

\subsection{SIR Model}

The SIR for $f$-FAMA is given by
\begin{equation} \label{eq:SIR_BC_fFAMA}
    \text{SIR}= \underset{n}{\text{max}} \frac{\sigma_s^2 \lvert h_n^{(u,u)} \rvert^2}
    { \lvert \widetilde h_n^{(u,u)} \rvert^2}
    = \underset{n}{\text{max}} \frac{X_n}{ \widehat U  Z_n}
\end{equation}
in which $\widehat U = \widetilde U / \widetilde m = \widetilde \sigma^2 / \sigma^2 \sigma_s^2$, $X_n$ is given by (\ref{eq:Xk}) and $Z_n$ is defined as
\begin{align} \label{eq:Zn}
    Z_n  =  
     \sum_{l=1}^{\widetilde m}  
     \left(\widetilde{x}^{(u)}_{nl}+ \varphi \, \widetilde{\mathrm{x}}^{(u)}_{b(n)l} \right)^2 
    + \left(\widetilde{y}^{(u)}_{nl}+ \varphi \, \widetilde{\mathrm{y}}^{(u)}_{b(n)l} \right)^2.
\end{align}

\textit{Remark 1:} 
Comparing the SIR for $f$-FAMA in (\ref{eq:SIR_BC_fFAMA}) and the SIR for $s$-FAMA in (\ref{eq:SIR_BC}), note that the expressions are similar, except for the constant $\widehat{U}$ and the upper limits of the sum of $Y_k$ and $Z_k$. 
As a consequence, the SIR-based OP for $f$-FAMA can be obtained from the SIR-based OP for $s$-FAMA in (\ref{eq:OP_SIR_sFAMA}) substituting $\gamma$ by $\widehat U \gamma$ and $\widetilde U$ by $\widetilde m$.

%\textit{Remark 3:} The multiplexing gain for $f$-FAMA is the same as that of $s$-FAMA, but considering the OP of $f$-FAMA.
\textit{Remark 2:} The multiplexing gain for $f$-FAMA is calculated as~\eqref{eq:mult_gain}, but considering the OP particularities of $f$-FAMA.

% \textit{Remark 3:} As $\lvert \widetilde h^{(u)}_{n} \rvert^2 = \sum_{l=1}^{\widetilde m} \lvert  \widetilde{g}^{(u)}_{nl} \rvert ^2$,  $\widetilde m$ is restricted to integers that depend on the values of $m_{\tilde u}$ and $U$ in (\ref{eq:mTil}). 
% \hl{However, the SIR-based OP for $f$-FAMA does not have mathematical restrictions for non-integer $\widetilde m$. }
% \textcolor{red}{Está errado, pois $\widetilde m$ é o limite de um somatório, então tem que ser inteiro. A expressão da SINR é que permite $\widetilde m$ não ser inteiro, pois está na Bessel}

\section{Opportunistic FAMA under Nakagami-$m$ Channels with Spatial Block-Correlation Model}
\label{Ofama}

In opportunistic scheduling, the best $U$ users are selected from a pool of $M (\geq U)$ users to maximize network capacity.
We combine in this section opportunistic scheduling with FAMA, operating in either fast or slow modes, referred to as O-FAMA. In FAMA, the BS performs no pre-processing, and a user’s OP is independent of the others.
Consequently, in O-FAMA, the BS can identify the top users by sequentially activating and deactivating them until the best set is found. In this paper, we assume that the strongest $U$ UEs are selected as an ideal condition. In~\cite{Waqar2024_Opportunistic_FAMA_RL}, it is shown that it is possible to select the best users and FAS ports to reach a network sum rate close to the ideal based in a reinforcement learning approach, which makes O-FAMA an interesting option for multiple access.

% In opportunistic scheduling, the best $U$ users are selected from a poll of $M(\geq U)$ users to maximize the network capacity.
% In this section, opportunistic scheduling is combined with FAMA, operating in fast or slow modes, in what is called O-FAMA.
% In FAMA, there is no pre-processing in the BS and the OP of one user does not depend on the others.
% Therefore, in O-FAMA the BS can select the top users simply turning them on and off until the best are found.
% In contrast, other types of scheduling shall test all possible combinations of users to find the network best capacity. 

% In this paper, we assume that the strongest $U$ UEs are selected as an ideal condition.
% Using a reinforcement learning approach, in \cite{Waqar2024_Opportunistic_FAMA_RL} it was shown that it is possible to select the best users and FAS ports to reach a network sum rate close to the ideal, which makes O-FAMA an interesting option for multiple access.

Based on the selection of the $U$ users with better channel conditions from a pool $M$ users, the multiplexing gain of an O-FAMA network is given by \cite[Eq. (21-a)]{Wong2023OpportunisticFAMA} $\mathcal{G}_m = \sum_{u=M-U+1}^M \mathcal{I}_{1-P_{\text{out}}} (M-u+1,u),$
% \begin{equation} \label{eq:mult_gain_OFAMA}
%     \mathcal{G}_m = \sum_{u=M-U+1}^M \mathcal{I}_{1-P_{\text{out}}} (M-u+1,u),
% \end{equation}
where $\mathcal{I}_x(a,b)$ is the regularized incomplete beta function. Note that $\mathcal{G}_m $ depends on the OP, so the multiplexing gain is dependent on whether it is operating in fast or slow mode. Moreover, $\mathcal{G}_m$ can be approximated by \cite[Eq. (46)]{Wong2023OpportunisticFAMA} $\mathcal{G}_m \approx \text{min} \{ U, M(1 - P_{\text {out}})   \}. $

\section{Numerical Results}
\label{results}

This section presents the numerical results for the performance metrics developed in this work. Monte Carlo simulations are also employed to validate the derived analytical expressions\footnote{The code is available at: https://github.com/HigoTh/famablc.}. The Gauss–Laguerre quadrature approximations are evaluated with 
$n_{I}=n_{J}=50$ roots. In addition, the block-correlation analysis is based on a reference correlation matrix derived from the Jakes model.
Our results compare the exact OP with the analytical approximation obtained via Gauss–Laguerre quadrature and the Monte Carlo simulations. 
All curves exhibit a strong overlap, validating our analysis.

\begin{figure}[!htb]
    \centering
    \includegraphics[width=0.9\linewidth]{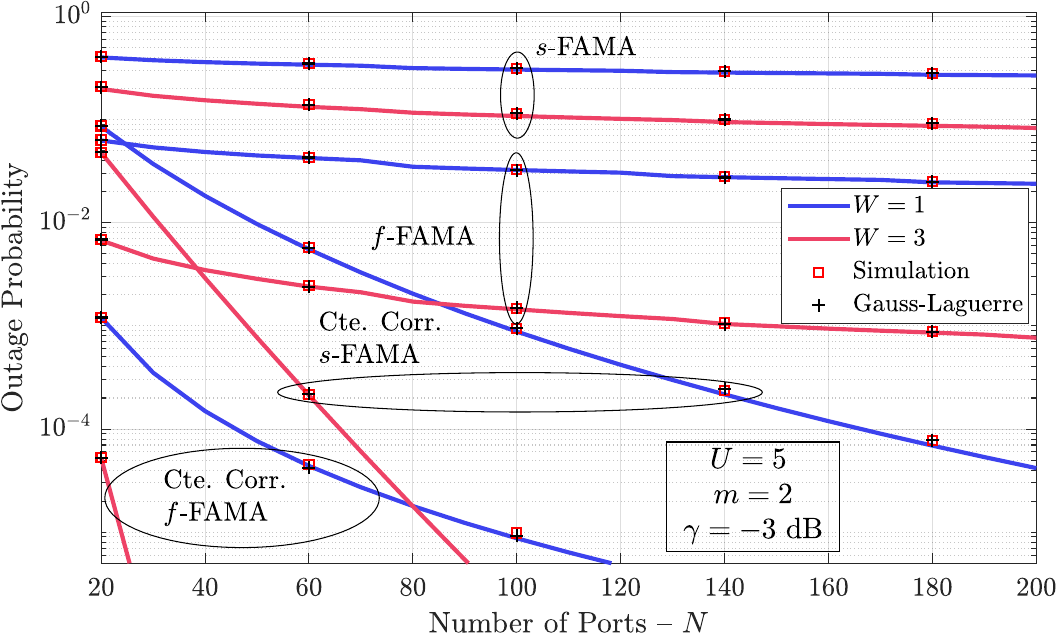}
    \caption{OP curves for $s$-FAMA and $f$-FAMA systems, considering different $W$ and correlation models, with $U=5$, $m = 2$, and $\gamma = -3$~dB.
    %\hl{Incluir elipse e limitante inferior da OP}
    }
    \label{fig:Fig1}
\end{figure}

Fig.~\ref{fig:Fig1} shows the OP curves as a function of the number of ports $N$ for both $s$-FAMA and $f$-FAMA systems, considering different values of the normalized antenna length $W$. The remaining parameters are fixed at $U=5$ users, $m=2$ and 
$\gamma=-3$~{dB}. This result highlights the impact of different correlation models on OP, comparing Jakes-based block correlation with constant correlation across ports, presented in~\cite{Wong2022Closed}. 
The results indicate that, for equivalent systems, FAMA achieves considerably better performance in fast mode, owing to port selection optimization at the symbol time scale. Increasing the number of ports also improves performance by reducing OP. However, this gain is significantly diminished under a more realistic block-correlation model. 
The constant-correlation model reproduces a non-realistic scenario with rapid signal fluctuations across ports, offering more opportunities for SIR maximization and thus leading to an overestimation of performance. Conversely, with the block-correlation model, spatial fluctuations are slower, reflecting the practical characteristics of spatial correlation, where fewer opportunities for SIR maximization occurs and consequently limiting the OP improvement achievable with larger $N$.
%\textcolor{red}{Paulo: In the constant correlation model, for a given $W$, low $N$ results in a large space between ports and indeed low correlation, but as $N$ increases, it also increases the correlation and the inaccuracy of the model.}
Finally, it is noted that for both the FAMA systems and the adopted correlation models, the OP improves with the increase of $W$ since a larger $W$ reduces the spatial correlation among ports.

% OP curves as a function of the number of ports $N$ for both $s$-FAMA and $f$-FAMA systems are illustrated in Fig.~\ref{fig:Fig1}, under $W=1$ and $W=3$, considering different correlation models, $U=5$, $m = 2$, and $\gamma = -3$~dB.
% Our results compare the exact OP with the analytical approximation obtained via Gauss–Laguerre quadrature and the Monte Carlo simulations. The three curves exhibit an almost perfect overlap, thereby confirming the tightness of the approximation and validating its accuracy against both the theoretical benchmark and empirical results.
% Note that OP reduces as the number of ports $N$ increases and for fixed $N$, better performance is obtained for $f$-FAMA under the two correlation models considered.
% For $s$-FAMA or $f$-FAMA and considering the block-correlation model, a plateau in the OP curves is perceived due to the saturation effect caused by the strong spatial correlation between closed ports.
% Fig.~\ref{fig:Fig1} also shows as a benchmark the OP considering a constant correlation model presented in~\cite{Wong2022Closed}, which underestimates the performance particularly for high values of $N$. This simplified correlation model misleadingly optimistic performance, while the block-correlation approach remains accurate and analytically tractable.
% Finally, it is noted that for both the FAMA systems and the adopted correlation models, the OP improves with the increase of $W$.
% This occurs because a larger $W$ reduces the spatial correlation among ports, thus improving the performance.

\begin{figure}[!htb]
    \centering
    \includegraphics[width=0.9\linewidth]{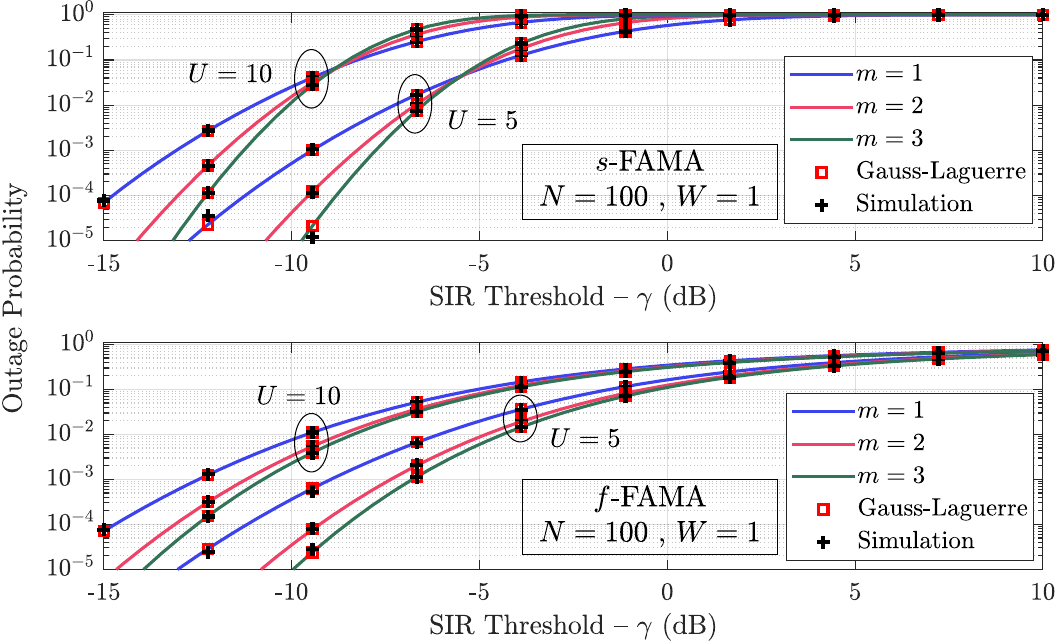}
    \caption{OP curves as a function of the SIR threshold for $s$-FAMA and $f$-FAMA considering different values of $m$ and $U$, with $N=100$ and $W=1$.
    %\hl{Centralizar o nome Outage Probability na figura e incluir o limitante superior no gráfico}
    }
    \label{fig:Fig2}
\end{figure}

Fig.~\ref{fig:Fig2} presents the OP curves as a function of the SIR threshold for $s$-FAMA and $f$-FAMA under different values of the fading parameter $m$ and the number of users $U$. The antenna parameters are fixed at $N=100$ and $W=1$. As a benchmark, the Rayleigh case ($m=1$) is included. As expected, the OP increases with higher SIR threshold requirements. Moreover, $f$-FAMA consistently outperforms $s$-FAMA under the same system configuration. The impact of $m$ differs across the evaluated cases. For $s$-FAMA, under low SIR threshold requirements, the Rayleigh case ($m=1$) represents the worst performance. However, at higher SIR thresholds, the case with $m=3$ becomes the worst. For $f$-FAMA, the OP exhibits a general improvement as $m$ increases.
Furthermore, $U$ worsens the OP since higher interference leads to higher outage levels.

% OP curves as a function of the SIR threshold are presented in Fig.~\ref{fig:Fig2} for $s$-FAMA and $f$-FAMA under different fading parameters $m$ and number of users $U$, considering $N=100$ and $W=1$.
% Note that a strong adherence is observed between the theoretical, simulated, and approximate results, which corroborates our mathematical analysis.
% As a benchmark, the Rayleigh case ($m=1$) is presented.
% In both scenarios, (i) when increasing the SIR threshold $\gamma$ the OP tends to 1; (ii) an increase in $m$ improves performance, given that the system's fading intensity is lower and (iii) increasing $U$ worsens the OP, since higher interference leads to higher outage levels.
% For the same OP, we perceive that $f$-FAMA supports a much larger number of users and tolerates higher SIR thresholds.

\begin{figure}[!htb]
    \centering
    \includegraphics[width=0.9\linewidth]{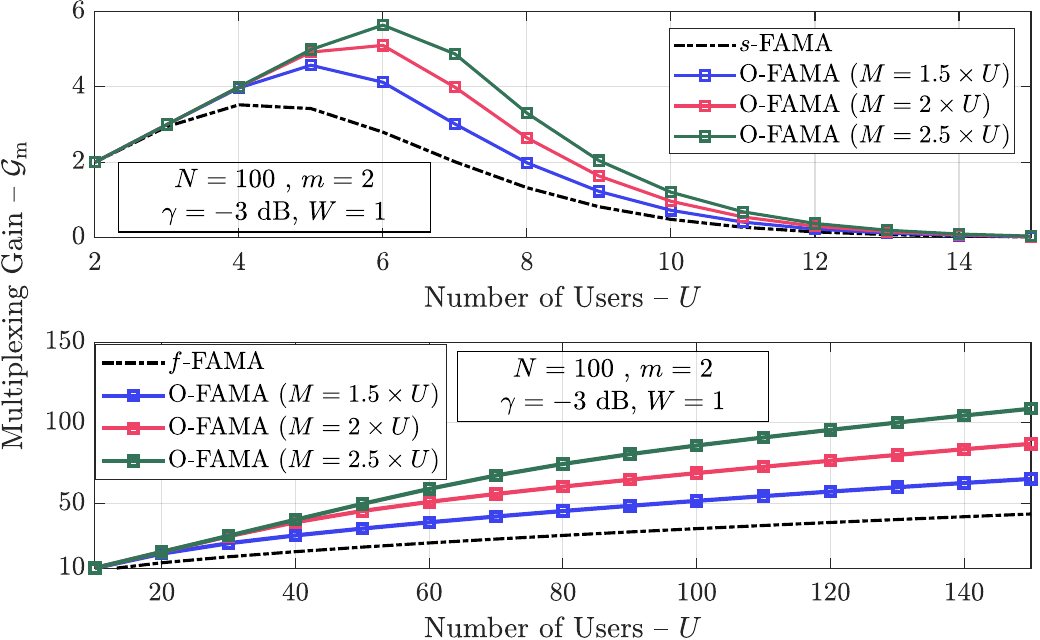}
    \caption{Multiplexing gain curves as a function of the number of users $U$ for $s$-FAMA, $f$-FAMA and O-FAMA under different settings of $M$, considering
$N=100$, $W=1$, $m = 2$, and $\gamma = -3$~dB.}
    \label{fig:Fig3}
\end{figure}

Fig.~\ref{fig:Fig3} depicts the multiplexing gain curves as a function of the number of users for $s$-FAMA, $f$-FAMA, and O-FAMA (in slow and fast modes). The system parameters are set to $N=100$ ports, $m=2$, $\gamma=-3$~{dB}, and $W=1$. The performance of O-FAMA is evaluated under different user pool sizes $M$. As observed, increasing $M$ while keeping $U$ fixed enhances the multiplexing gain, indicating that a larger number of users can be served without experiencing outage. Conversely, for FAMA in slow mode, an inflection point emerges, beyond which further increases in network size reduce the multiplexing gain. In contrast, in fast mode this effect does not occur, and the multiplexing gain exhibits a consistently increasing trend with network scaling. This notable performance advantage is attributed to the ability to maximize the SIR at symbol time, thereby offering more optimization opportunities compared to the slow mode.

% Fig.~\ref{fig:Fig3} presents multiplexing gain curves as a function of the number of users $U$ for $s$-FAMA, $f$-FAMA and O-FAMA under different settings of $M$, considering
% $N=100$, $W=1$, $m = 2$, and $\gamma = -3$~dB.
% The results reveal that O-FAMA consistently outperforms $s$-FAMA and $f$-FAMA across the entire range of users, particularly when $M$ is scaled proportionally with $U$.
% For the upper subfigure, note that the gain peaks around $U = 6$ users, reaching its maximum value at $M = 2.5 \times U$. After this peak, $\mathcal{G}_m$ gradually decreases as the number of users increases, indicating the impact of multiuser interference.
% The bottom subfigure highlights the improvement of $\mathcal{G}_m$ with $U$.
% Notably, the performance gap between O-FAMA and
% $f$-FAMA also widens as $U$ increases, underscoring the scalability and robustness of the O-FAMA scheme.
% For $s$-FAMA and $f$-FAMA, note from~(\ref{eq:mult_gain}) that the maximum $\mathcal{G}_m$ is $U$. For $s$-FAMA the maximum is obtained with a small number of users, whereas in $f$-FAMA the system still does not reach this value even for $U = 100$ users.

\section{Conclusion}
\label{conclusions}

This paper studied slow, fast, and opportunistic FAMA under Nakagami-$m$ fading channels, considering the spatial block-correlation model, where OP expressions and bounds for slow and fast FAMA were derived and applied to the study O-FAMA.
We obtained multiplexing gain results for O-FAMA over fast and slow FAMA. 
Numerical and simulation results confirmed the accuracy of the proposed analytical framework.
Our work extends FAMA analysis to more general fading and correlation conditions, providing practical insights.
%into the design of interference-limited multiple access systems. 

%Our work extends FAMA analysis to more general fading and correlation conditions but also provides practical insights into the design of interference-limited multiple access systems. 
%The results demonstrate the superiority of opportunistic scheduling, especially when combined with fast FAMA, making it a promising candidate for future wireless networks.

\bibliography{refs_IEEEstyleV2}
\bibliographystyle{IEEEtran}

\end{document}